\documentclass[12pt]{amsart}
\usepackage[margin=2.0cm]{geometry}

\usepackage{amsfonts}
\usepackage{mathrsfs}
\usepackage{cite}
\usepackage{enumerate}
\usepackage{xcolor}

\newcommand{\R}{{\mathbb R}}

\newtheorem{theorem}{Theorem}[section]
\newtheorem{lemma}[theorem]{Lemma}
\newtheorem{prop}[theorem]{Proposition}
\newtheorem{corollary}[theorem]{Corollary}

\newcommand{\sfe}{S^{n-1}}

\title[The Orlicz version of the $L_p$ Minkowski problem for $-n<p<0$]{The Orlicz version of the $L_p$ Minkowski problem  
\\
for $-n<p<0$}
\author[Gabriele Bianchi, K\'aroly J. B\"or\"oczky, Andrea Colesanti]
{Gabriele Bianchi, K\'aroly J. B\"or\"oczky and Andrea Colesanti}
\address{Dipartimento di Matematica e Informatica ``U. Dini", Universit\`a di Firenze, Viale Morgagni 67/A, Firenze, Italy I-50134} \email{gabriele.bianchi@unifi.it}
\address{Alfr\'ed R\'enyi Institute of Mathematics, Hungarian Academy
  of Sciences, Reltanoda u. 13-15, H-1053 Budapest, Hungary, and
Department of Mathematics, Central European University, Nador u 9, H-1051, Budapest, Hungary} \email{boroczky.karoly.j@renyi.mta.hu}
\address{Dipartimento di Matematica e Informatica ``U. Dini", Universit\`a di Firenze, Viale Morgagni 67/A, Firenze, Italy I-50134} \email{andrea.colesanti@unifi.it}
\thanks{First and third authors are supported in part by the Gruppo Nazionale per l'Analisi Matematica, la Probabilit\`a e le loro Applicazioni (GNAMPA) of the Istituto Nazionale di Alta Matematica (INdAM). Second author is supported in part by
NKFIH grants 116451, 121649 and 129630.}

\begin{document}

\begin{abstract}
Given a function $f$ on the unit sphere $S^{n-1}$, the $L_p$ Minkowski problem asks for a convex body $K$ whose $L_p$ surface area measure has density $f$ with respect to the standard $(n-1)$-Hausdorff measure on $S^{n-1}$. 
In this paper we deal with the generalization of this problem which arises in the Orlicz-Brunn-Minkowski theory when an Orlicz function $\varphi$ substitutes the $L_p$ norm and  $p$ is  in the range $(-n,0)$.
This problem is equivalent to solve the Monge-Ampere equation
\[
 \varphi(h)\det(\nabla^2 h+h I)=f\quad\text{ on $S^{n-1}$,}
\]
where $h$ is the support function of the convex body $K$.
\end{abstract}

\maketitle

\section{Introduction}

We work in the $n$-dimensional Euclidean space $\R^n$, $n\ge2$. A \emph{convex body} $K$ in $\mathbb{R}^n$ is a 
compact convex set that has non-empty interior. Given a convex body $K$, for $x\in\partial K$ we denote by $\nu_K(x)\subset S^{n-1}$ the family of all unit exterior normal vectors to $K$ at $x$  (the {\em Gau{\ss} map}). We can then define the \emph{surface area measure} $S_{K}$ of $K$, which is a Borel measure on the unit sphere $S^{n-1}$ of ${\mathbb R}^n$, as follows: for a Borel set $\omega\subset S^{n-1}$ we set
$$
S_{K}(\omega)=\mathcal{H}^{n-1}\left(\nu_K^{-1}(\omega)\right)=
\mathcal{H}^{n-1}\left(\{x\in\partial K:\,\nu_K(x)\cap\omega\neq\emptyset\}\right)
$$
(see, {\em e.g.}, Schneider \cite{SCH}). 

The classical Minkowski problem can be formulated as follows: {\em given a Borel measure $\mu$ on $\sfe$, find a convex  body $K$ such that $\mu=S_K$}. The reader is referred to \cite[Chapter 8]{SCH} for an exhaustive presentation of this problem and its solution. 

Throughout this paper we will consider (either for the classical Minkowski problem or for its variants) the case in which $\mu$ has a density $f$ with respect to the $(n-1)$-dimensional Hausdorff measure on $\sfe$. Under this assumption the Minkowski problem is equivalent to solve (in the classic or in the weak sense) a differential equation on the sphere. Namely:
\begin{equation}
\label{intro 1}
\det(\nabla^2h+h I)=f,
\end{equation}
where: $h$ is the support function of $K$, $\nabla^2h$ is the matrix formed by the second covariant derivatives of $h$ with respect to a local orthonormal frame on $\sfe$ and $I$ is the identity matrix of order $(n-1)$. 

Many different types of variations of the Minkowski problem have been considered (we refer for instance to \cite[Chapters 8 and 9]{SCH}). Of particular interest for our purposes is the so called $L_p$ version of the problem. At the origin of this new problem there is the replacement of the usual Minkowski addition of convex bodies by the $p$-addition. As an effect, the corresponding differential equation takes the form
\begin{equation}
\label{intro 2}
h^{1-p}\det(\nabla^2h+h I)=f,
\end{equation}
(see \cite[Section 9.2]{SCH}). The study of the $L_p$ Minkowski problems developed in a significant way in the last decades, as a part of the so called $L_p$ Brunn-Minkowski theory, which represents now a substantial area of Convex Geometry. One of the most interesting aspects of this problem is that several threshold values of the parameter $p$ can be identified, {\em e.g.} $p=1$, $p=0$, $p=-n$, across which the nature of the problem changes drastically. For an account on the literature and on the state of the art of 
the $L_p$ Minkowski problem (especially for the values $p<1$) we refer the reader to \cite{BBCY} and \cite{BiBoCo}. 

Of particular interest here is the range $-n<p<0$. In this case Chou and Wang (see \cite{CW}) solved the corresponding problem when the measure $\mu$ has a density $f$, and $f$ is bounded and bounded away from zero. This result was slightly generalised by the authors in collaboration with Yang in \cite{BBCY}, where $f$ is allowed to be in $L_{\frac{n}{n+p}}$. 

\begin{theorem}[Chou and Wang; Bianchi, B\"or\"oczky, Colesanti and Yang]
\label{0pn}
For $n\geq 1$ and $-n<p<0$, if the non-negative and non-trivial function $f$ is in $L_{\frac{n}{n+p}}(S^{n-1})$ then
\eqref{intro 2} has a solution in the Alexandrov sense; namely, $f\,d\mathcal{H}^{n-1}=dS_{K,p}$ for a convex body 
$K\in \mathcal{K}_0^n$. In addition, if $f$ is invariant under a closed subgroup $G$ of $O(n)$, then $K$ can be chosen to be invariant under $G$.
\end{theorem}

As a further extension of the $L_p$ Minkowski problem, one may consider its Orlicz version. Formally, this problem arises in the context of the Orlicz-Brunn-Minkowski theory of convex bodies (see \cite[Chapter 9]{SCH}). In practice, the relevant differential equation is 
\begin{equation*}
\varphi(h)\det(\nabla^2 h+h I)=f,
\end{equation*}
where $\varphi$ is a suitable {\em Orlicz function}. The $L_p$ Minkowski problem is obtained when $\varphi(t)=t^{1-p}$, for $t\ge0$. 

When $\varphi\colon(0,\infty)\to(0,\infty)$ is continuous and monotone decreasing, this problem (under a symmetry assumption) has been considered by Haberl, Lutwak, Yang, Zhang in \cite{HLYZ10}. Comparing the previous assumptions on $\varphi$ with the $L_p$ case, we see that this corresponds to the values $p\ge 1$. 

We are interested in the case in which the monotonicity assumption is reversed, corresponding to the values $p<1$. Hence we assume that $\varphi\colon[0,\infty)\to\R$ is continuous and monotone increasing, having the example 
$\varphi(t)=t^{1-p}$, $p<1$, as a prototype. To control in a more precise form the behaviour of $\varphi$ with that of a power function, we assume that there exists $p<1$ such that
\begin{equation}\label{intro 3}
\liminf_{t\to 0^+}\frac{\varphi(t)}{t^{1-p}}>0.
\end{equation}
Concerning the behaviour of $\varphi$ at $\infty$ we impose the condition:
\begin{equation}\label{intro 4}
\int_1^\infty\frac{1}{\varphi(t)}\,dt<\infty.
\end{equation}

The corresponding Minkowski problem in this setting can be called the Orlicz $L_p$ Minkowski problem. The solution of this problem in the range $p\in(0,1)$ is due to Jian, Lu \cite{JLZ18+}. We also note that Orlicz versions of the so called $L_p$ dual Minkowski have been considered recently by Gardner, Hug, Weil, Xing, Ye \cite{GHWXY}, Gardner, Hug, Xing, Ye \cite{GHXY}, Xing, Ye, Zhu \cite{XYZ} and Xing, Ye \cite{XY}. 

In this paper we focus on the range of values $p\in(-n,0)$. As an extension of the results contained in \cite{BBCY}, we establish the following existence theorem (note that, as usual in the case of Orlicz versions of Minkowski type problems, we can only provide a solution up to a constant factor). 

\begin{theorem}
\label{psi0pn}
For $n\geq 2$, $-n<p<0$ and monotone increasing continuous function $\varphi:[0,\infty)\to[0,\infty)$
satisfying $\varphi(0)=0$, and conditions \eqref{intro 3} and \eqref{intro 4}, if the non-negative non-trivial 
function $f$ is in $L_{\frac{n}{n+p}}(S^{n-1})$, then there exists $\lambda>0$ and a convex body 
$K\in \mathcal{K}_{0}^n$ with $V(K)=1$ such that
$$
\lambda\varphi(h)\det(\nabla^2 h+h I)=f
$$
holds for $h=h_K$ in the Alexandrov sense; namely, 
$\lambda\varphi(h_K)\,dS_{K}=f\,d\mathcal{H}^{n-1}$. In addition, if $f$ is invariant under a closed subgroup $G$ of $O(n)$, then $K$ can be chosen to be invariant under $G$.
\end{theorem}

We note that the origin may lie on $\partial K$ for the solution $K$ in Theorem~\ref{psi0pn}.

We observe that Theorem~\ref{psi0pn} readily yields Theorem~\ref{0pn}. Indeed if 
$-n<p<0$, $f\in L_{\frac{n}{n+p}}(S^{n-1})$, $f\geq0$, $f\not\equiv0$ and 
$\lambda h_K^{1-p}\,dS_{K}=f\,d\mathcal{H}^{n-1}$ for
$K\in \mathcal{K}_{0}^n$ and $\lambda>0$, then 
$h_{\widetilde{K}}^{1-p}\,dS_{\widetilde{K}}=f\,d\mathcal{H}^{n-1}$
for $\widetilde{K}=\lambda^{\frac1{n-p}}K$.

In Section~\ref{secstrategy} we sketch the proof of Theorem~\ref{psi0pn} and describe the structure of the paper.

\section{Notation}

The scalar product on $\R^n$ is denoted by $\langle\cdot,\cdot\rangle$, and the corresponding Euclidean norm is denoted by $\|\cdot \|$. The $k$-dimensional Hausdorff measure normalized in such a way that it coincides with the Lebesgue measure on $\R^k$ is denoted by $\mathcal{H}^{k}$. The angle (spherical distance) of $u,v\in S^{n-1}$ is denoted by $\angle(u,v)$.

We write $\mathcal{K}_0^n$ ($\mathcal{K}_{(0)}^n$)  to denote the family of convex bodies with $o\in K$ ($o\in{\rm int}\,K$). Given a convex body $K$, for a Borel set $\omega\subset S^{n-1}$, $\nu_K^{-1}(\omega)$ is the Borel set of $x\in \partial K$ with $\nu_K(x)\cap \omega\neq \emptyset$. A point $x\in\partial K$ is called smooth if $\nu_K(x)$ consists of a unique vector, and in this case, we use $\nu_K(x)$ to denote this unique vector, as well. It is well-known that $\mathcal{H}^{n-1}$-almost every $x\in\partial K$ is smooth (see, {\em e.g.}, Schneider \cite{SCH}); let $\partial'K$ denote the family of smooth points of $\partial K$.

For a convex compact set $K$ in $\R^n$, let $h_K$ be its support 
function:
$$
h_K(u)=\max\{\langle x,u\rangle:\, x\in K\} \mbox{ \ \ for $u\in\R^n$}.
$$
Note that if $K\in{\mathcal K}_0^n$, then $h_K\ge0$. If $p\in\R$ and 
$K\in {\mathcal K}_0^n$, then the  $L_{p}$-surface area measure is defined by
$$
dS_{K,p}=h_K^{1-p}\,d S_K
$$
where for $p>1$ the right-hand side is assumed to be a finite measure.
In particular, if $p=1$, then $S_{K,p}=S_K$, and if $p<1$ and $\omega\subset S^{n-1}$ Borel, then
$$
S_{K,p}(\omega)=\int_{\nu_{K}^{-1}(\omega)}\langle x,\nu_{K}(x)\rangle^{1-p}d\mathcal{H}^{n-1}(x).
$$

\section{Sketch of the proof of Theorem~\ref{psi0pn}}
\label{secstrategy}

To sketch the argument leading to Theorem~\ref{psi0pn}, first we consider
the case when $-n<p\leq-(n-1)$ and $\varphi(t)=t^{1-p}$, and $\tau_1\leq f\leq \tau_2$ 
for some constants $\tau_2>\tau_1>0$. 
We set $\psi(t)=1/\varphi(t)=t^{p-1}$ for $t>0$,
and define $\Psi:(0,\infty)\to(0,\infty)$   by
$$
\Psi(t)=\int_t^\infty \psi(s)\,ds=\mbox{$-\frac{1}{p}$}\,t^p,
$$
which is a strictly convex function.

Given a convex body $K$ in $\R^n$, we set
$$
\Phi(K,\xi)=\int_{S^{n-1}}\Psi(h_{K-\xi})f\,d\mathcal{H}^{n-1};
$$
this is a strictly convex function of $\xi\in{\rm int}\,K$. 
As $f>\tau_1$ and $p\leq -(n-1)$, there is a (unique) $\xi(K)\in {\rm int}\,K$ such that
$$
\Phi(K,\xi(K))=\min_{\xi\in {\rm int}\,K}\Phi(K,\xi).
$$
This statement is proved in Proposition~\ref{xiinside}, but the conditions $f>\tau_1$ and $p\leq -(n-1)$
are actually used in the preparatory statement Lemma~\ref{insidegood}.

Using $p>-n$ and the Blaschke-Santal\'o inequality
(see 
Lemma~\ref{extremal}
and the preparatory statement Lemma~\ref{diameterest}), one verifies that
there exists a convex body $K_0$ in $\R^n$ with $V(K_0)=1$ maximizing 
$\Phi(K,\xi(K))$ over all convex bodies $K$ in $\R^n$ with $V(K)=1$.

Finally a variational argument proves that there exists $\lambda_0>0$
such that $f\,d\mathcal{H}^{n-1}=\lambda_0\varphi(h_{K_0})\,dS_{K_0}$.
A crucial ingredient (see Lemma~\ref{xider}) is that, as $\psi$ is $C^1$ and $\psi'<0$,
$\Phi(K_t,\xi(K_t))$ is a differentiable function of $K_t$ for a suitable variation $K_t$ of $K_0$.

In the general case, when still keeping the condition $\tau_1\leq f\leq \tau_2$
but allowing any $\varphi$ which satisfies the assumptions of Theorem~\ref{psi0pn},
we meet two main obstacles.
On the one hand, even if $\varphi(t)=t^{1-p}$ but $0<t<-(n-1)$, it may happen that
for a convex body $K$ in $\R^n$, the infimum of $\Phi(K,\xi)$ for $\xi\in{\rm int}\,K$ is attained when 
$\xi$ tends to the boundary of $K$. On the other hand, the possible lack of differentiability of $\varphi$ (or equivalently of $\psi$) destroys the variational argument.

Therefore, we approximate $\psi$ by smooth functions, and also make sure that the
approximating functions are large enough near zero
 to ensure that the minimum of the analogues of $\Phi(K,\xi)$ as a function of $\xi\in {\rm int}\,K$ exists
for any convex body $K$. 

Section~\ref{secpreliminary} proves some preparatory statements, Section~\ref{secenergy} introduces 
the suitable analogue of the energy function $\Phi(K,\xi(K))$, and Section~\ref{secvariational} provides the variational formula  for an extremal body for the energy function.
We prove Theorem~\ref{psi0pn} if $f$ is bounded and bounded away from zero in Section~\ref{secmaintheorembounded}, and finally in full strength in Section~\ref{secmaintheorem}.

\section{Some preliminary estimates}
\label{secpreliminary}

In this section, we prove the simple but technical estimates 
Lemmas~\ref{bigpsibehave} and \ref{diameterest}
that will be used in various settings during the main argument.

\begin{lemma}
\label{bigpsibehave}
For $\delta\in(0,1)$, $A,\tilde{\aleph}>0$ and $q\in(-n,0)$,
let $\widetilde{\psi}:\,(0,\infty)\to(0,\infty)$ satisfy that
$\widetilde{\psi}(t)\leq \tilde{\aleph} t^{q-1}$ for $t\in(0,\delta]$ and
$\int_\delta^\infty \widetilde{\psi}\leq A$. If $t\in(0,\delta)$ and
$\tilde{\aleph}_0=\max\{\frac{\tilde{\aleph}}{|q|},\frac{A}{\delta^q}\}$, then
$\widetilde{\Psi}(t)=\int_t^\infty \widetilde{\psi}$ satisfies
$$
\widetilde{\Psi}(t)\leq\tilde{\aleph}_0 t^q.
$$
\end{lemma}
\proof We observe that if $t\in(0,\delta)$, then
$$
\widetilde{\Psi}(t)\leq \int_t^{\delta} \tilde{\psi}(s)\,ds+A\leq
\tilde{\aleph}\int_t^{\delta} s^{q-1}\,ds+A=
\frac{\tilde{\aleph}}{|q|}(t^q-\delta^q)+A\leq
t^q \max\left\{\frac{\tilde{\aleph}}{|q|},\frac{A}{\delta^q}\right\}.
\mbox{ \ \ Q.E.D.}
$$

We write $B^n$ to denote the Euclidean unit ball in $\R^n$, and set $\kappa_n=\mathcal{H}^n(B^n)$.
For a convex body $K$ in $\R^n$, let $\sigma(K)$ denote its centroid, which satisfies (see Schneider \cite{SCH})
\begin{equation}
\label{centroid}
-(K-\sigma(K))\subset n(K-\sigma(K)).
\end{equation}

Next, if $o\in {\rm int}\,K$ then the polar of $K$ is
$$
K=\{x\in\R^n:\,\langle x,y\rangle\leq 1\;\forall y\in K\}=
\{tu:\,u\in S^{n-1}\mbox{ and }0\leq t\leq h_K(u)^{-1}\}.
$$
In particular, the Blaschke-Santal\'o inequality $V(K) V((K-\sigma(K))^*)\leq V(B^n)^2$ 
(see Schneider \cite{SCH})  yields that
\begin{equation}
\label{BS}
\int_{S^{n-1}}h_{K-\sigma(K)}^{-n}d\mathcal{H}^{n-1}\leq \frac{nV(B^n)^2}{V(K)}.
\end{equation}

As a preparation for the proof of Lemma~\ref{diameterest}, we need the following statement about absolutely continuous measures. For $t\in(0,1)$ and $v\in S^{n-1}$, we consider the spherical strip
$$
\Xi(v,t)=\{u\in S^{n-1}\,:\,|\langle u,v\rangle|\leq t\}.
$$

\begin{lemma}
\label{strips}
If $f\in L_1(S^{n-1})$ and
$$
\varrho_f(t)=\sup_{v\in S^{n-1}}\int_{\Xi(v,t)}|f|\,d\mathcal{H}^{n-1}
$$
for $t\in(0,1)$, then we have $\lim_{t\to 0^+}\varrho_f(t)=0$.
\end{lemma}
\proof We observe that $\varrho_f(t)$ is decreasing, therefore the limit 
$\lim_{t\to 0^+}\varrho_f(t)=\delta\geq 0$ exists.
We suppose that $\delta>0$, and seek a contradiction.

Let $\mu$ be the absolutely continuous measure $d\mu=\,|f|\,d\mathcal{H}^{n-1}$ on $S^{n-1}$.
According to the definition of $\varrho_f$,  for any $k\geq 2$, there exists some $v_k\in S^{n-1}$ such that
$\mu(\Xi(v_k,\frac1k)) \geq \delta/2$,
Let $v\in S^{n-1}$ be an accumulation point of the sequence $\{v_k\}$. For any $m\geq 2$, there exists
$\alpha_m>0$ such that $\Xi(u,\frac1{2m})\subset \Xi(v,\frac1m)$ if $u\in S^{n-1}$ and $\angle(u,v)\leq \alpha_m$.
Since for any $m\geq 2$, there exists some $k\geq 2m$ such that
$\angle(v_k,v)\leq \alpha_m$,  we have
$\mu(\Xi(v,\frac1m))\geq \mu(\Xi(v_k,\frac1k))\geq \delta/2$. We deduce that
$\mu(v^\bot\cap S^{n-1})=\mu\left(\cap_{m\geq 2}\Xi(v,\frac1m)\right)\geq \delta/2$, which contradicts
$\mu(v^\bot\cap S^{n-1})=0$.
\hfill Q.E.D.

 \begin{lemma}
\label{diameterest}
For $\delta\in(0,1)$,  $\tilde{\aleph}>0$ and  $q\in(-n,0)$, let 
$\widetilde{\Psi}:\,(0,\infty)\to(0,\infty)$ be a monotone 
decreasing continuous function such that 
$\widetilde{\Psi}(t)\leq \tilde{\aleph} t^q$ for $t\in(0,\delta]$ and $\lim_{t\to\infty}\widetilde{\Psi}(t)=0$, and let 
$\tilde{f}$ be a non-negative function in $L_{\frac{n}{n+p}}(S^{n-1})$.
Then for any $\zeta>0$, there exists 
a $D_\zeta$ depending on $\zeta$,
$\widetilde{\Psi}$, $\delta$, $\tilde{\aleph}$, $q$ and  $\tilde{f}$
such that
if $K$ is a convex body in $\R^n$ with 
$V(K)=1$ and ${\rm diam}\,K\geq D_\zeta$ then
$$
\int_{S^{n-1}}(\widetilde{\Psi}\circ h_{K-\sigma(K)})\,\tilde{f}\,d\mathcal{H}^{n-1}\leq \zeta.
$$
\end{lemma}
\proof We may assume that $\sigma(K)=o$. Let  $R=\max_{x\in K}\|x\|$, and let $v\in S^{n-1}$ such that
$Rv\in K$. It follows from (\ref{centroid}) that $-\frac{R}{n}\,v\in K$. 

Since $\lim_{t\to\infty}\widetilde{\Psi}(t)=0$ and $\tilde{f}$ is in $L_1(S^{n-1})$ by the H\"older inequality, we can choose $r\geq 1$ such that
\begin{equation}
\label{rcond}
\widetilde{\Psi}(r)\int_{S^{n-1}} \tilde{f}\,d\mathcal{H}^{n-1}<\frac{\zeta}2.
\end{equation}
We partition $S^{n-1}$ into the two measurable parts
\begin{eqnarray*}
\Xi_0&=&\{u\in S^{n-1}:\,h_K(u)\geq r\}\\
\Xi_1&=&\{u\in S^{n-1}:\,h_K(u)< r\}.
\end{eqnarray*} 

Let us  estimate the integrals over $\Xi_0$ and $\Xi_1$. We deduce from (\ref{rcond}) that
\begin{equation}
\label{Xi0}
\int_{\Xi_0}(\widetilde{\Psi}\circ h_{K})\,\tilde{f}\,d\mathcal{H}^{n-1}\leq \frac{\zeta}2.
\end{equation}

Next we claim that
\begin{equation}
\label{Xi1in}
\Xi_1\subset \Xi\left(v,\frac{nr}{R}\right).
\end{equation}
For any $u\in \Xi_1$, we choose $\eta\in\{-1,1\}$ such that $\langle u,\eta v\rangle\geq 0$,
thus $\frac{\eta R}n\, v\in K$ yields that 
$r>h_K(u)\geq \langle u,\frac{\eta R}n\, v\rangle$. In turn, we conclude (\ref{Xi1in}).
It follows from (\ref{Xi1in}) and Lemma~\ref{strips} that for the $L_1$ function $f=\tilde{f}^{\frac{n}{n+q}}$, we have
\begin{equation}
\label{Xi1f}
\int_{\Xi_1}\tilde{f}^{\frac{n}{n+q}} \leq \varrho_f\left(\frac{nr}{R}\right).
\end{equation}

To estimate the decreasing function $\widetilde{\Psi}$ on $(0,r)$, we claim that if $t\in (0,r)$ then
\begin{equation}
\label{bigpsion0r}
\widetilde{\Psi}(t)\leq \frac{\tilde{\aleph} \delta^q}{r^q}\,t^q.
\end{equation}
We recall that $r\geq 1>\delta$. In particular, if $t\leq \delta$, then $\widetilde{\Psi}(t)\leq\tilde{\aleph} t^q$ yields (\ref{bigpsion0r}).
If $t\in(\delta,r)$, then using that $\widetilde{\Psi}$ is decreasing, (\ref{bigpsion0r}) follows from
$$
\widetilde{\Psi}(t)\leq \widetilde{\Psi}(\delta)\leq\frac{\tilde{\aleph} \delta^q}{t^q}\,t^q
\leq \frac{\tilde{\aleph} \delta^q}{r^q}\,t^q.
$$

Applying first (\ref{bigpsion0r}), then the H\"older inequality, after that
the Blaschke-Santal\'o inequality (\ref{BS}) with $V(K)=1$ and finally (\ref{Xi1f}),
we deduce that
\begin{eqnarray*}
\int_{\Xi_1}(\widetilde{\Psi}\circ h_{K})\,\tilde{f}\,d\mathcal{H}^{n-1}&\leq&
\frac{\tilde{\aleph} \delta^q}{r^q}
\int_{\Xi_1}h_{K}^{-|q|}\,\tilde{f}\,d\mathcal{H}^{n-1}\\
&\leq&
\frac{\tilde{\aleph} \delta^q}{r^q}\left(\int_{\Xi_1}h_{K}^{-n}\,d\mathcal{H}^{n-1}\right)^{\frac{|q|}{n}}
\left(\int_{\Xi_1}\tilde{f}^{\frac{n}{n-|q|}}\,d\mathcal{H}^{n-1}\right)^{\frac{n-|q|}{n}}\\
&\leq&
\frac{\tilde{\aleph} \delta^q}{r^q}\left(nV(B^n)^2\right)^{\frac{|q|}{n}}
\varrho_f\left(\frac{nr}{R}\right)^{\frac{n+q}{n}}.
\end{eqnarray*}
Therefore after fixing $r\geq 1$ satisfying (\ref{rcond}), we may choose $R_0> r$ 
such that
$$
\frac{\tilde{\aleph} \delta^q}{r^q}\,n^{\frac{|q|}{n}}V(B^n)^{\frac{2|q|}{n}}
\varrho_f\left(\frac{nr}{R_0}\right)^{\frac{n+q}{n}}
<\frac{\zeta}2
$$
by Lemma~\ref{diameterest}. In particular, if $R\geq R_0$, then
$$
\int_{\Xi_1}(\widetilde{\Psi}\circ h_{K})\,\tilde{f}\,d\mathcal{H}^{n-1}\leq \frac{\zeta}2.
$$
Combining this estimate with (\ref{Xi0}) shows that setting $D_\zeta=2R_0$, if 
${\rm diam}\,K\geq D_\zeta$, then $R\geq R_0$, and hence
$\int_{S^{n-1}}(\widetilde{\Psi}\circ h_{K})\,\tilde{f}\,d\mathcal{H}^{n-1}\leq \zeta$.
\hfill Q.E.D. 

\section{The extremal problem related to Theorem~\ref{psi0pn} when  $f$ is bounded and bounded away from zero}
\label{secenergy}

For $0<\tau_1<\tau_2$,  let the real function $f$ on $S^{n-1}$ satisfy
 \begin{equation}
\label{ftau12}
\tau_1<f(u)<\tau_2\mbox{ \ for $u\in S^{n-1}$}.
\end{equation}
In addition, let $\varphi:[0,\infty)\to[0,\infty)$ be a 
 continuous monotone increasing  function satisfying $\varphi(0)=0$, 
$$
\liminf_{t\to 0^+}\frac{\varphi(t)}{t^{1-p}}>0
\mbox{ \ and \ }\int_1^\infty\frac{1}{\varphi(t)}\,dt<\infty.
$$
It will be more convenient to work with the decreasing function $\psi=1/\varphi:(0,\infty)\to(0,\infty)$, which has the properties
\begin{eqnarray}
\label{psicond0}
\limsup_{t\to 0^+}\frac{\psi(t)}{t^{p-1}}&<&\infty\\
\label{psicondinf}
\int_1^\infty\psi(t)\,dt&<&\infty.
\end{eqnarray}
We consider the function $\Psi:(0,\infty)\to(0,\infty)$  defined by
$$
\Psi(t)=\int_t^\infty \psi(s)\,ds,
$$
which readily satisfies
\begin{eqnarray}
\label{bigpsider}
\Psi'&=&-\psi,\mbox{ and hence $\Psi$ is convex and strictly monotone decreasing},\\
\label{bigpsiinfty}
\lim_{t\to \infty}\Psi(t)&=&0.
\end{eqnarray}
According to (\ref{psicond0}),
there exist some $\delta\in(0,1)$ and $\aleph>1$ such that
\begin{equation}
\label{psicondaleph1}
\psi(t)< \aleph t^{p-1} \mbox{ \ for $t\in(0,\delta)$}.
\end{equation}

As we pointed out in Section~\ref{secstrategy}, we smoothen $\psi$ using convolution.
Let   $\eta:\,\R \to [0,\infty)$ be a non-negative $C^\infty$ ``approximation of identity''
with ${\rm supp}\,\eta\subset[-1,0]$
and $\int_\R\eta=1$. For any $\varepsilon\in(0,1)$, we consider the non-negative  
$\eta_\varepsilon(t)=
\frac1{\varepsilon}\eta(\frac{t}{\varepsilon})$
satisfying that $\int_\R\eta_\varepsilon=1$, ${\rm supp}\,\eta_\varepsilon\subset[-\varepsilon,0]$, and define
$\theta_\varepsilon:(0,\infty)\to(0,\infty)$ by
$$
\theta_\varepsilon(t)=\int_\R \psi(t-\tau)\eta_\varepsilon(\tau)\,d\tau=\int_{-\varepsilon}^0 \psi(t-\tau)\eta_\varepsilon(\tau)\,d\tau.
$$
As $\psi$ is monotone decreasing and continuous on $(0,\infty)$, the properties of $\eta_\varepsilon$ yield
\begin{eqnarray*}
\theta_\varepsilon(t)&\leq &\psi(t)\mbox{ \ for $t>0$ and $\varepsilon\in(0,1)$}\\
\theta_\varepsilon(t_1)&\geq &\theta_\varepsilon(t_2) \mbox{ \ for $t_2>t_1>0$ and $\varepsilon\in(0,1)$}\\
\theta_\varepsilon&\mbox{tends}&\mbox{uniformly to $\psi$ on any interval with positive endpoints as $\varepsilon$ tends to zero.}
\end{eqnarray*}
Next, for any $t_0>0$, the function $l_{t_0}$ on $\R$ defined by
$$
l_{t_0}(t)=
\left\{\begin{array}{rcl}
\psi(t)&\mbox{ if }&t\geq t_0\\
0&\mbox{ if }& t<t_0
\end{array}\right.
$$
is bounded, and  hence locally integrable. For the convolution $l_{t_0} * \eta_\varepsilon$, we have that
 $(l_{t_0} * \eta_\varepsilon)(t)=\theta_\varepsilon(t)$ for $t>t_0$ and $\varepsilon\in(0,1)$, thus
$$
\mbox{$\theta_\varepsilon$ is $C^1$ for each $\varepsilon\in(0,1)$.}
$$

As it is explained in Section~\ref{secstrategy}, we need to modify $\psi$ in a way such that the new function
is of order at least $t^{-(n-1)}$ if $t>0$ is small. We set
$$
q=\min\{p,-(n-1)\},
$$
and hence (\ref{psicondaleph1}) and $\delta\in(0,1)$ yields that
\begin{equation}
\label{psicondaleph}
\theta_\varepsilon(t)\leq \psi(t)< \aleph t^{q-1} \mbox{ \ for $t\in(0,\delta)$ and $\varepsilon\in(0,\delta)$}.
\end{equation}

Next we construct $\tilde{\theta}_\varepsilon:(0,\infty)\to(0,\infty)$ satisfying
\begin{eqnarray*}
\tilde{\theta}_\varepsilon(t)&=&\theta_\varepsilon(t)\leq\psi(t)\mbox{ \ for $t\geq \varepsilon$ and $\varepsilon\in(0,\delta)$}\\
\tilde{\theta}_\varepsilon(t)&\leq &\aleph t^{q-1}\mbox{ \ for $t\in(0,\delta)$ and 
$\varepsilon\in(0,\delta)$}\\
\tilde{\theta}_\varepsilon(t)&=&\aleph t^{q-1}\mbox{ \ for $t\in(0,\frac{\varepsilon}2]$ and 
$\varepsilon\in(0,\delta)$}\\
\tilde{\theta}_\varepsilon&\mbox{ is $C^1$ }&\mbox{ and is monotone decreasing.}
\end{eqnarray*}
It follows that  
$$
\mbox{$\tilde{\theta}_\varepsilon$
tends uniformly to $\psi$ on any interval with positive endpoints as $\varepsilon$ tends to zero.}
$$
To construct suitable $\tilde{\theta}_\varepsilon$,
first we observe that the conditions above determine $\tilde{\theta}_\varepsilon$ outside
the interval $(\frac{\varepsilon}2,\varepsilon)$, and
$\tilde{\theta}_\varepsilon(\varepsilon)<\aleph\varepsilon^{q-1}$. Writing $\Delta$ to denote the degree one polynomial whose graph is the tangent to the graph of $t\mapsto \aleph t^{q-1}$ at $t=\varepsilon/2$, we have
$\Delta(t)<\aleph t^{q-1}$ for $t>\varepsilon/2$ and $\Delta(\varepsilon)<0$. Therefore we can choose
$t_0\in(\frac{\varepsilon}2,\varepsilon)$ such that 
$\tilde{\theta}_\varepsilon(\varepsilon)<\Delta(t_0)<\aleph\varepsilon^{q-1}$.
We define $\tilde{\theta}_\varepsilon(t)=\Delta(t)$ for $t\in (\frac{\varepsilon}2,t_0)$, and construct 
$\tilde{\theta}_\varepsilon$ on $(t_0,\varepsilon)$ in a way that $\tilde{\theta}_\varepsilon$ stays $C^1$
on $(0,\infty)$. It follows from the way $\tilde{\theta}_\varepsilon$ is constructed that
$\tilde{\theta}_\varepsilon(t)\leq \aleph t^{q-1}$ also for $t\in [\frac{\varepsilon}2,\varepsilon]$.

In order to ensure a negative derivative, we consider
$\psi_\varepsilon:(0,\infty)\to(0,\infty)$ defined by
\begin{equation}
\psi_\varepsilon(t)=\tilde{\theta}_\varepsilon(t)+\frac{\varepsilon}{1+t^2}
\end{equation}
for $\varepsilon\in(0,\delta)$ and $t>0$.
This $C^1$ function $\psi_\varepsilon$ has the following properties:
\begin{equation}
\label{psiepsilon}
\begin{array}{rcll}
\psi_\varepsilon(t)&\leq &\psi(t)+\frac1{1+t^2}&\mbox{for $t\geq \varepsilon$ and $\varepsilon\in(0,\delta)$}\\[0.5ex]
\psi'_\varepsilon(t)&< &0& \mbox{for $t>0$ and $\varepsilon\in(0,\delta)$}\\[0.5ex]
\psi_\varepsilon(t)&< &2\aleph t^{q-1}&\mbox{for $t\in(0,\delta)$ and 
$\varepsilon\in(0,\delta)$}\\[0.5ex]
\psi_\varepsilon(t)&> &\aleph t^{q-1}&\mbox{for $t\in(0,\frac{\varepsilon}2)$ and 
$\varepsilon\in(0,\delta)$}\\[0.5ex]
\psi_\varepsilon&\mbox{tends}&\mbox{uniformly}&\mbox{to $\psi$ on any interval with positive endpoints as 
$\varepsilon$ tends to zero.}
\end{array}
\end{equation}

For $\varepsilon\in(0,\delta)$, we also consider the $C^2$ function $\Psi_\varepsilon:(0,\infty)\to(0,\infty)$  defined by
$$
\Psi_\varepsilon(t)=\int_t^\infty \psi_\varepsilon(s)\,ds,
$$
and hence (\ref{psiepsilon})  yields
\begin{eqnarray}
\label{bigpsieinfty}
\lim_{t\to \infty}\Psi_\varepsilon(t)&=&0\\
\label{bigpsieder}
\Psi'_\varepsilon&=&-\psi_\varepsilon,\mbox{ \ thus $\Psi_\varepsilon$ is strictly decreasing and strictly convex}.
\end{eqnarray}

For $\varepsilon\in(0,\delta)$, 
Lemma~\ref{bigpsibehave} and
(\ref{psiepsilon}) imply that setting
$$
A=\int_{\delta}^\infty \psi(t)+\frac1{1+t^2}\,dt, 
$$
we have
\begin{equation}
\label{bigpsicondup}
\Psi_\varepsilon(t)
\leq\aleph_0 t^q\mbox{ \ for $\aleph_0=\max\{\frac{2\aleph}{|q|},\frac{A}{\delta^q}\}$ and $t\in(0,\delta)$}.
\end{equation}
On the other hand, if $\varepsilon\in(0,\delta)$ and $t\in (0,\frac{\varepsilon}4)$, then
\begin{equation}
\label{bigpsicondlow}
\Psi_\varepsilon(t)\geq \int_t^{\varepsilon/2}\aleph s^{q-1}\,ds=
\frac{\aleph}{|q|}(t^q-(\varepsilon/2)^q)\geq \frac{\aleph}{|q|}(t^q-(2t)^q)=
\aleph_1 t^q\mbox{ \ for $\aleph_1=\frac{(1-2^q)\aleph}{|q|}>0$}.
\end{equation}

According to (\ref{psiepsilon}), we have $\lim_{\varepsilon\to 0^+}\psi_\varepsilon(t)=\psi(t)$ 
and $\psi_\varepsilon(t)\leq \psi(t)+\frac1{1+t^2}$ for
any $t>0$, therefore Lebesgue's Dominated Convergence Theorem implies
\begin{equation}
\label{bigpsilimit}
\lim_{\varepsilon\to 0^+}\Psi_\varepsilon(t)=\Psi(t) 
\mbox{ \  for any $t>0$}.
\end{equation}
It also follows from (\ref{psiepsilon}) that if $t\geq \varepsilon$, then
\begin{equation}
\label{bigpsiupper}
\Psi_\varepsilon(t)=\int_t^\infty\psi_\varepsilon\leq \int_t^\infty \psi(s)+\frac1{1+s^2}\,ds
\le\Psi(t) +\frac{\pi}2.
\end{equation}

For any convex body $K$ and $\xi\in{\rm int}\,K$, we consider
$$
\Phi_\varepsilon(K,\xi)=\int_{S^{n-1}}(\Psi_\varepsilon\circ h_{K-\xi})f\,d\mathcal{H}^{n-1}=
\int_{S^{n-1}}\Psi_\varepsilon(h_{K}(u)-\langle \xi,u\rangle)f(u)\,d\mathcal{H}^{n-1}(u).
$$
Naturally, $\Phi_\varepsilon(K)$ depends on $\psi$ and $f$, as well, but we do not signal these dependences.

We equip $\mathcal{K}_{0}^n$ with the Hausdorff metric, which is the $C_\infty$ metric on the space of the restrictions of support functions to $S^{n-1}$. For $v\in S^{n-1}$ and $\alpha\in[0,\frac{\pi}2]$, we consider the spherical cap
$$
\Omega(v,\alpha)=\{u\in S^{n-1}\,\langle u,v\rangle\geq \cos\alpha\}.
$$
We write $\pi:\R^n\backslash \{o\}\to S^{n-1}$ the radial projection:
$$
\pi(x)=\frac x{\|x\|}.
$$
In particular, if $\pi$ is restricted to the boundary of a $K\in \mathcal{K}_{(0)}^n$, then this map is Lipschitz. Another typical application of the radial projection
is to consider, for $v\in S^{n-1}$, the composition $x\mapsto \pi(x+v)$  as a map $v^\bot\to S^{n-1}$  where  
\begin{equation}
\label{Jacobian}
\mbox{the Jacobian of $x\mapsto \pi(x+v)$ at $x\in v^\bot$ is $(1+\|x\|^2)^{-n/2}$.}
\end{equation}

The following Lemma~\ref{insidegood} is the statement where we apply directly that $\psi$ is modified to be essentially $t^q$ if $t$ is very small.

\begin{lemma}
\label{insidegood}
Let $\varepsilon\in(0,\delta)$, and let $\{K_i\}$ be a sequence of convex bodies tending to a convex body $K$ in $\R^n$, and let
$\xi_i\in{\rm int}\,K_i$ such that $\lim_{i\to \infty}\xi_i=x_0\in \partial K$. Then
$$
\lim_{i\to \infty}\Phi_\varepsilon(K_i,\xi_i)=\infty.
$$
\end{lemma}
\proof We may assume that $\lim_{i\to \infty}\xi_i=x_0=o$.  
Let $v\in S^{n-1}$ be an exterior normal to  $\partial K$ at $o$, and choose some $R>0$ such that $K\subset RB^n$. Therefore we may assume that $K_i-\xi_i\subset (R+1)B^n$, $h_{K_i}(v)<\varepsilon/8$ 
and $\|\xi_i\|<\varepsilon/8$ for all $\xi_i$, thus $h_{K_i-\xi_i}(v)<\varepsilon/4$ for all $i$.

For any $\zeta\in(0,\frac{\varepsilon}8)$, there exists $I_\zeta$ such that if $i\geq I_\zeta$, then $\|\xi_i\|\leq \zeta/2$ and $\langle y,v\rangle\leq \zeta/2$ for all $y\in K_i$,
and hence $\langle y,v\rangle\leq \zeta$ for all $y\in K_i-\xi_i$.
 For $i\geq I_\zeta$, any  $y\in K_i-\xi_i$ can be written in the form $y=sv+z$ where $s\leq \zeta$ and 
$z\in v^\bot\cap(R+1)B^n$, thus if $\angle(v,u)=\alpha\in[\zeta,\frac{\pi}2)$ for $u\in S^{n-1}$, then we have
\begin{equation}
\label{hKualpha}
h_{K_i-\xi_i}(u)\leq (R+1)\sin\alpha+\zeta\cos\alpha\leq (R+2)\alpha.
\end{equation}
We set $\beta=\frac{\varepsilon}{4(R+2)}$, and for $\zeta\in (0,\beta)$, we define
$$
\Omega_\zeta=\Omega(v,\beta)\backslash\Omega(v,\zeta).
$$
In particular, as $\Psi_\varepsilon(t)\geq \aleph_1 t^q$ for $t\in (0,\frac{\varepsilon}4)$ according to
(\ref{bigpsicondlow}), if $u\in \Omega_\zeta$, then (\ref{hKualpha}) implies
$$
\Psi_\varepsilon(h_{K_i-\xi_i}(u))\geq \gamma (\angle(v,u))^q
$$
for $\gamma=\aleph_1(R+2)^q$.

The function $x\mapsto \pi(x+v)$ maps 
$B_\zeta =v^\bot\cap \big((\tan\beta) B^n\backslash (\tan\zeta) B^n\big)$ bijectively onto 
$\Omega_\zeta$, while $\beta<\frac18$ and  (\ref{Jacobian}) yield that the
 Jacobian of this map is at least $2^{-n}$ on $B_\zeta$. 

Since $f>\tau_1$ and $\angle(v,\pi(x+v))\leq 2x$ for $x\in B_\zeta$, if $i\geq I_\zeta$, then
\begin{eqnarray*}
\Phi_\varepsilon(K_i,\xi_i)&=&\int_{S^{n-1}}\Psi_\varepsilon(h_{K_i-\xi_i}(u))f(u)\,d\mathcal{H}^{n-1}(u)
\geq \int_{\Omega_\zeta}\tau_1\gamma  (\angle(v,u))^q\,d\mathcal{H}^{n-1}(u)\\
&\geq&
 \frac{\tau_1\gamma}{2^{n+|q|}}\int_{B_\zeta}\|x\|^q\,d\mathcal{H}^{n-1}(x)=
\frac{(n-1)\kappa_{n-1}\tau_1\gamma}{2^{n+|q|}}\int_{\tan\zeta}^{\tan \beta}t^{q+n-2}\,dt.
\end{eqnarray*}
As $\zeta>0$ is arbitrarily small and $q\leq 1-n$, we conclude that
$\lim_{i\to\infty}\Phi_\varepsilon(K_i,\xi_i)=\infty$. \hfill Q.E.D. \\

Now we single out the optimal $\xi\in{\rm int}\,K$.

\begin{prop}
\label{xiinside}
For $\varepsilon\in(0,\delta)$ and a convex body $K$ in $\R^n$, there exists a unique $\xi(K)\in{\rm int}\,K$ such that
$$
\Phi_\varepsilon(K,\xi(K))=\min_{\xi\in{\rm int}\,K}\Phi_\varepsilon(K,\xi).
$$
In addition, $\xi(K)$ and $\Phi_\varepsilon(K,\xi(K))$ are continuous functions of $K$, and 
$\Phi_\varepsilon(K,\xi(K))$ is translation invariant.
\end{prop}
\proof 
The first part of this proof, the one regarding the existence of $\xi(K)\in{\rm int}\,K$ and its uniqueness, is very similar to the proof of \cite[Proposition 3.2]{BBCY} given by the authors and Yang for  the $L_p$ Minkowski problem. It is very short and we rewrite it here for completeness. 

Let $\xi_1,\xi_2\in{\rm int}\,K$, $\xi_1\neq\xi_2$, and let $\lambda\in(0,1)$. 
If $u\in S^{n-1}\backslash(\xi_1-\xi_2)^\bot$, then $\langle u,\xi_1\rangle\neq\langle u,\xi_2\rangle$, and hence
the  strict convexity of $\Psi_\varepsilon$ (see (\ref{bigpsieder})) yields that
$$
\Psi_\varepsilon(h_K(u)-\langle u,\lambda \xi_1+(1-\lambda)\xi_2\rangle)>
\lambda\Psi_\varepsilon(h_K(u)-\langle u,\xi_1\rangle)+
(1-\lambda)\Psi_\varepsilon(h_K(u)-\langle u,\xi_2\rangle),
$$
thus $\Phi_\varepsilon(K,\xi)$ is a strictly convex function of $\xi\in{\rm int}\,K$ by $f>\tau_1$.

Let $\xi_i\in{\rm int}\,K$ such that 
$$
\lim_{i\to\infty}\Phi_\varepsilon(K,\xi_i)=\inf_{\xi\in{\rm int}\,K}\Phi_\varepsilon(K,\xi).
$$
We may assume that $\lim_{i\to\infty}\xi_i=x_0\in K$, and Lemma~\ref{insidegood} yields
 $x_0\in{\rm int}\,K$.
Since $\Phi_\varepsilon(K,\xi)$ is a strictly convex and continuous function of $\xi\in{\rm int}\,K$, 
$x_0$ is the unique minimum point of $\xi\mapsto \Phi_\varepsilon(K,\xi)$, which we denote by $\xi(K)$ (not signalling the dependence on $\varepsilon$, $\psi$ and $f$).

Readily $\xi(K)$ is translation equivariant, and $\Phi_\varepsilon(K,\xi(K))$ is translation invariant.

For the continuity of $\xi(K)$ and $\Phi_\varepsilon(K,\xi(K))$, let us
consider a sequence $\{K_i\}$ of convex bodies tending to a convex body $K$ in $\R^n$.
We may assume that $\xi(K_i)$ tends to a $x_0\in K$.

For any $y\in {\rm int}\, K$, there exists an $I$ such that $y\in{\rm int}\,K_i$ for $i\geq I$.
Since $h_{K_i}$ tends uniformly to $h_K$ on $S^{n-1}$, we have that
$$
\limsup_{i\to\infty}\Phi_\varepsilon(K_i,\xi(K_i))\leq \lim_{i\to\infty\atop i\geq I} 
\Phi_\varepsilon(K_i,y)=\Phi_\varepsilon(K,y).
$$
Again Lemma~\ref{insidegood} implies that $x_0\in {\rm int}\, K$. It follows that 
$h_{K_i-\xi_i(K_i)}$ tends uniformly to $h_{K-x_0}$, thus
$$
\Phi_\varepsilon(K,x_0)=\lim_{i\to\infty}\Phi_\varepsilon(K_i,\xi(K_i))\leq \lim_{i\to\infty\atop i\geq I} 
\Phi_\varepsilon(K_i,y)=\Phi_\varepsilon(K,y).
$$ 
In particular, $\Phi_\varepsilon(K,x_0)\leq \Phi_\varepsilon(K,y)$ for any $y\in {\rm int}\, K$, thus
$x_0=\xi(K)$. In turn, we deduce  $\xi(K_i)$ tends to $\xi(K)$, and $\Phi_\varepsilon(K_i,\xi(K_i))$
tends to $\Phi_\varepsilon(K,\xi(K))$.
Q.E.D.\\

Since $\xi\mapsto \Phi_\varepsilon(K,\xi)=
\int_{S^{n-1}}\Psi_\varepsilon(h_K(u)-\langle u,\xi\rangle)f(u)\,d\mathcal{H}^{n-1}(u)$ is maximal at 
$\xi(K)\in{\rm int}\,K$ and $\Psi'_\varepsilon=-\psi_\varepsilon$, we deduce

\begin{corollary}
\label{intcond}
For $\varepsilon\in(0,\delta)$ and a convex body $K$ in $\R^n$, we have
$$
\int_{S^{n-1}}u\ \psi_\varepsilon\Big(h_{K}(u)-\langle u,\xi(K)\rangle\Big)
f(u)\,d\mathcal{H}^{n-1}(u)=o.
$$
\end{corollary}

For a  closed subgroup $G$ of $O(n)$, we write $\mathcal{K}_{(0)}^{n,G}$ to denote the family of 
$K\in \mathcal{K}_{(0)}^n$ invariant under $G$.

\begin{lemma}
\label{extremal}
For $\varepsilon\in(0,\delta)$,
there exists a $K^\varepsilon\in \mathcal{K}_{(0)}^n$ with $V(K^\varepsilon)=1$ such that
$$
\Phi_\varepsilon(K^\varepsilon,\xi(K^\varepsilon))\geq \Phi_\varepsilon(K,\xi(K))
\mbox{ \ for any $K\in \mathcal{K}_{(0)}^n$ with $V(K)=1$}.
$$
In addition, if $f$ is invariant under a closed subgroup $G$ of $O(n)$, then $K^\varepsilon$ can be chosen to be invariant under $G$.
\end{lemma}
\proof We choose a sequence $K_i\in \mathcal{K}_{(0)}^n$ with $V(K_i)=1$ for $i\geq 1$ such that
$$
\lim_{i\to \infty}\Phi(K_i,\xi(K_i))=\sup\{\Phi(K,\xi(K)):\mbox{$K\in \mathcal{K}_{(0)}^n$ with $V(K)=1$}\}.
$$
Writing $B_1=\kappa_n^{-1/n}B^n$ to denote the unit ball centred at the origin and having volume $1$, we may assume that each $K_i$ satisfies
\begin{equation}
\label{Betabound}
\Phi_\varepsilon(K_i,\sigma(K_i))\geq\Phi_\varepsilon(K_i,\xi(K_i))\geq \Phi_\varepsilon(B_1,\xi(B_1)).
\end{equation}
According to Proposition~\ref{xiinside}, we may also assume that $\sigma_{K_i}=o$ for each $K_i$.

We deduce from Lemma~\ref{diameterest}, (\ref{bigpsieinfty}), (\ref{bigpsicondup}) and (\ref{Betabound}) that there exists some $R>0$ such that
$K_i\subset RB^n$ for any $i\geq 1$. According to the Blaschke selection theorem, we may assume that
$K_i$ tends to a compact convex set $K^\varepsilon$ with $o\in K^\varepsilon$. It follows from the continuity of the volume that $V(K^\varepsilon)=1$, and hence ${\rm int}\,K^\varepsilon\neq \emptyset$. We conclude from Lemma~\ref{xiinside}
 that $\Phi_\varepsilon(K^\varepsilon,\xi(K^\varepsilon))=\lim_{i\to \infty}\Phi_\varepsilon(K_i,\xi(K_i))$. 

If $f$ is invariant under a closed subgroup $G$ of $O(n)$, then we apply the same argument to convex bodies in
$\mathcal{K}_{(0)}^{n,G}$ instead of $\mathcal{K}_{(0)}^n$.
\hfill Q.E.D. \\

Since $\Phi(5)<\Phi(4)$, (\ref{bigpsilimit}) yields some $\tilde{\delta}\in(0,\delta)$ such that
$\Psi_\varepsilon(4)\geq \Phi(5)$ for $\varepsilon\in(0,\tilde{\delta})$.
For future reference, the monotonicity of 
$\Psi_\varepsilon$, ${\rm diam}\kappa_n^{-1/n}B^n\leq 4$ and (\ref{Betabound}) yield that
if $\varepsilon\in(0,\tilde{\delta})$, then
\begin{equation}
\label{control}
\Phi_\varepsilon(K^\varepsilon,\sigma(K^\varepsilon))\geq
\Phi_\varepsilon\big(\kappa_n^{-1/n}B^n,\xi(\kappa_n^{-1/n}B^n)\big)\geq
\int_{S^{n-1}}\Psi_\varepsilon(4)f\,d\mathcal{H}^{n-1}\geq
\Psi(5)\int_{S^{n-1}}f\,d\mathcal{H}^{n-1}.
\end{equation}

\section{Variational formulae and smoothness of the extremal body
 when $f$ is bounded and bounded away from zero}
\label{secvariational}

In this section, again let $0<\tau_1<\tau_2$ and  let the real function $f$ on $S^{n-1}$ satisfy
$\tau_1<f<\tau_2$.
In addition, let $\varphi$ be the continuous function  of Theorem~\ref{psi0pn}, and we use the notation developed in Section~\ref{secenergy}, say 
$\psi:(0,\infty)\to(0,\infty)$ is defined by
$\psi=1/\varphi$.

Now that we have constructed an extremal body $K^\varepsilon$, we want to show that it satisfies
the required differential equation in the Alexandrov sense by
using a variational argument. This section provides the formulae that  we will need, and ensure the required smoothness of $K^\varepsilon$.

Concerning the variation of volume, a key tool is
Alexandrov's Lemma~\ref{Alexandrov} (see Lemma 7.5.3 in \cite{SCH}). To state this,
for any continuous $h:S^{n-1}\to(0,\infty)$, we define the Alexandrov body
$$
[h]=\{x\in \R^n:\,\langle x,u\rangle \leq h(u)\mbox{ \ for }u\in S^{n-1}\}
$$
which is a convex body containing the origin in its interior.
Obviously, if $K\in \mathcal{K}_{(0)}^n$ then $K=[h_K]$.

\begin{lemma}[Alexandrov]
\label{Alexandrov}
For $K\in \mathcal{K}_{(0)}^n$ and a continuous function $g:\,S^{n-1}\to \R$,
$K(t)=[h_K+tg]$ satisfies
$$
\lim_{t\to 0}\frac{V(K(t))-V(K)}t=\int_{S^{n-1}}g(u)\,d S_{K}(u).
$$
\end{lemma}

To handle the variation of $\Phi_\varepsilon(K(t),\xi(K(t)))$ for a family $K(t)$ is a more subtle problem.
The next lemma shows essentially that if we perturb a convex body $K$ in a way such that the support function 
is  differentiable as a function of the parameter $t$ for $\mathcal{H}^{n-1}$-almost all $u\in S^{n-1}$, then
$\xi(K)$ changes also in a differentiable way. Lemma~\ref{xider} is the point of the proof where we use that
$\psi_\varepsilon$ is $C^1$ and $\psi'_\varepsilon<0$.

\begin{lemma} 
\label{xider}
For $\varepsilon\in(0,\delta)$, let $c>0$ and $t_0>0$, and let $K(t)$ be a family of convex bodies  with support function $h_t$ for $t\in[0,t_0)$. Assume that
\begin{description}
 \item[(i)] $|h_t(u)-h_0(u)|\leq c t$ for each $u\in S^{n-1}$ and $t\in[0,t_0)$,
 \item[(ii)] $\lim_{t\to 0^+}\frac{h_t(u)-h_0(u)}{t}$ exists for $\mathcal{H}^{n-1}$-almost all $u\in S^{n-1}$.
\end{description}
Then $\lim_{t\to 0^+}\frac{\xi(K(t))-\xi(K(0))}{t}$ exists.
 \end{lemma}
\proof We set $K=K(0)$. We may assume that $\xi(K)=o$, and hence Proposition~\ref{xiinside}
yields that
$$
\lim_{t\to 0^+}\xi(K(t))=o.
$$

There exists some $R>r>0$ such that $r\leq h_t(u)-\langle u,\xi(K(t))\rangle=h_{K(t)-\xi(K(t))}(u)\leq R$
for $u\in S^{n-1}$ and $t\in[0,t_0)$. Since $\psi_\varepsilon$ is $C^1$ on $[r,R]$, we can write
$$
\psi_\varepsilon(t)-\psi_\varepsilon(s)=\psi'_\varepsilon(s)(t-s)+\eta_0(s,t)(t-s)
$$
for $t,s\in[r,R]$ where $\lim_{t\to s}\eta_0(s,t)=0$.
Let $g(t,u)=h_t(u)-h_K(u)$ for $u\in S^{n-1}$ and $t\in[0,t_0)$.
Since $h_{K(t)-\xi(K(t))}$ tends uniformly to $h_K$ on $S^{n-1}$, we deduce that
if $t\in[0,t_0)$, then
\begin{equation}
\label{etudef}
\psi_\varepsilon\Big(h_t(u)-\langle u,\xi(K(t))\rangle\Big)-\psi_\varepsilon(h_K(u))=
\psi'_\varepsilon(h_K(u)) \big(g(t,u)-\langle u,\xi(K(t))\rangle\big)
+e(t,u)
\end{equation}
where
$$
|e(t,u)|\leq \eta(t)|g(t,u)-\langle u,\xi(K_t)\rangle|\mbox{ \ and \ }
\eta(t)=\eta_0(h_K(u),h_t(u)-\langle u,\xi(K(t))\rangle).
$$
Note that $\lim_{t\to 0^+}\eta(t)=0$ uniformly in $u\in S^{n-1}$.

In particular, (i) yields that $e(t,u)=e_1(t,u)+e_2(t,u)$ where 
\begin{equation}
\label{e1e2}
|e_1(t,u)|\leq c\eta(t)t\mbox{ \ and \ }|e_2(t,u)|\leq \eta(t)\|\xi(K(t))\|.
\end{equation}

It follows from (\ref{etudef}) and from applying Corollary~\ref{intcond} to $K(t)$ and $K$ that
$$
\int_{S^{n-1}}u\  \Big(\psi'_\varepsilon(h_K(u))\ \big(g(t,u)-\langle u,\xi(K(t))\rangle\big)
+e(t,u)\Big)\,f(u)d\mathcal{H}^{n-1}(u)=o,
$$
which can be written as
\begin{eqnarray*}
\int_{S^{n-1}}u\,  \psi'_\varepsilon(h_K(u))\,  g(t,u)\,f(u)d\mathcal{H}^{n-1}(u)&&\\
+\int_{S^{n-1}}u\,e_1(t,u)\,f(u)d\mathcal{H}^{n-1}(u)&=&
\int_{S^{n-1}}u\  \langle u,\xi(K_t)\rangle\psi'_\varepsilon(h_K(u))\,f(u)d\mathcal{H}^{n-1}(u)\\
&&-\int_{S^{n-1}}u\  e_2(t,u)\,f(u)d\mathcal{H}^{n-1}(u).
\end{eqnarray*}
Since $\psi'_\varepsilon(s)<0$ for all $s>0$, the symmetric matrix
$$
A=\int_{S^{n-1}}u\otimes u\, \psi'_\varepsilon(h_K(u))\,f(u)d\mathcal{H}^{n-1}(u)
$$
is negative definite because for any $v\in S^{n-1}$, we have
$$
v^TAv=\int_{S^{n-1}}\langle u,v\rangle^2\,\psi'_\varepsilon(h_K(u))\, f(u)\,d\mathcal{H}^{n-1}(u)<0.
$$
In addition, $A$ satisfies
$$
\int_{S^{n-1}}u\  \langle u,\xi(K_t)\rangle\,\psi'_\varepsilon(h_K(u))\,f(u)d\mathcal{H}^{n-1}(u)
=A \,\xi(K_t).
$$
It follows from (\ref{e1e2}) that if $t\geq 0$ is small, then
\begin{equation}
\label{psi1psi2}
A^{-1}\int_{S^{n-1}}u\  \psi'_\varepsilon(h_K(u))\  g(t,u)\,f(u)d\mathcal{H}^{n-1}(u)
+\tilde{e}_1(t)=\xi(K_t)-\tilde{e}_2(t),
\end{equation}
where $\|\tilde{e}_1(t)\|\leq \alpha_1 \eta(t) t$ and $\|\tilde{e}_2(t)\|\leq \alpha_2 \eta(t)\|\xi(K_t)\|$
for constants $\alpha_1,\alpha_2>0$. Since $\eta(t)$ tends to $0$ with $t$, 
if $t\geq 0$ is small, then $\|\xi(K(t))-\tilde{e}_2(t)\|\geq \frac12\,\|\xi(K_t)\|$.
Adding the estimate $g(t,u)\leq ct$, we deduce that
$\|\xi(K(t))\|\leq \beta\,t$ for a constant $\beta>0$, which in turn yields that
$\lim_{t\to 0^+}\frac{\|\tilde{e}_i(t)\|}t=0$ and $\tilde{e}_i(0)=0$ for $i=1,2$.
 Since
there exists $\lim_{t\to 0^+}\frac{g(t,u)-g(0,u)}{t}=\partial_1 g(0,u)$
for $\mathcal{H}^{n-1}$ almost all 
$u\in S^{n-1}$, and $\frac{g(t,u)-g(0,u)}{t}<c$ for all $u\in S^{n-1}$ and $t>0$, we conclude that
$$
\left.\frac{d}{dt}\,\xi(K(t))\right|_{t=0^+}=
A^{-1}\int_{S^{n-1}}u\  \psi'_\varepsilon(h_K(u))\  \partial_1g(0,u)\, f(u)\,d\mathcal{H}^{n-1}(u).
\mbox{ \ \ \ \ Q.E.D.}
$$

\begin{corollary}
\label{center-irrelevant}
Under the conditions of Lemma~\ref{xider}, and setting $K=K(0)$, we have
$$
\left.\frac{d}{dt}\,\Phi_\varepsilon(K(t),\xi(K(t)))\right|_{t=0^+}=
-\int_{S^{n-1}}\left.\frac{\partial}{\partial t} h_{K(t)}(u)\right|_{t=0^+} 
\psi_\varepsilon\big(h_K(u)-\langle u,\xi(K)\rangle\big)\,f(u)\,d\mathcal{H}^{n-1}(u).
$$
\end{corollary}
We omit the proof of this result since it  is very similar to that of \cite[Corollary 3.6]{BBCY}, given by the authors and Yang for  the $L_p$ Minkowski problem, with $f(u)\,d\mathcal{H}^{n-1}(u)$, $\Psi_\varepsilon$, $-\psi_\varepsilon$, Lemma~\ref{xider} and Corollary~\ref{intcond}  replacing respectively $d\mu(u)$, $\varphi_\varepsilon$, $\varphi'_\varepsilon$, Lemma 3.5 and Corollary 3.3. 
\bigskip

Given a family $K(t)$ of convex bodies for $t\in[0,t_0)$, $t_0>0$, 
to handle the variation of $\Phi_\varepsilon(K(t),\xi(K(t)))$ at $K(0)=K$ via applying
Corollary~\ref{center-irrelevant}, we need the properties (see Lemma~\ref{xider})
that there exists $c>0$ such that
\begin{eqnarray}
\label{xider1}
|h_{K(t)}(u)-h_K(u)|\leq c\,|t| &&\mbox{ \ for any $u\in S^{n-1}$ and $t\in[0,t_0)$}\\
\label{xider2}
\lim_{t\to 0^+}
\frac{h_{K(t)}(u)-h_K(u)}{t}&&\mbox{ \ exists for $\mathcal{H}^{n-1}$ almost all $u\in S^{n-1}$}.
\end{eqnarray}
 However, even if $K(t)=[h_K+th_C]$ for $K,C\in \mathcal{K}_{(0)}^n$ and for $t\in(-t_0,t_0)$,
 $K$ must satisfy some smoothness assumption in order to ensure that 
(\ref{xider2}) holds also for the two sided limits (problems occur say if $K$ is a polytope and $C$ is smooth). 

We recall that $\partial' K$  denotes the set of smooth points of $\partial K$.
We say that $K$ is \emph{quasi-smooth} if $\mathcal{H}^{n-1}(S^{n-1}\backslash\nu_K(\partial' K))=0$; 
namely, the set of  $u\in S^{n-1}$ that are exterior normals only at singular points has $\mathcal{H}^{n-1}$-measure zero.
The following Lemma~\ref{Wulffvariation}, taken from Bianchi, B\"or\"oczky, Colesanti, Yang \cite{BBCY}, shows that 
(\ref{xider1})
and (\ref{xider2}) are satisfied even if $t\in(-t_0,t_0)$ at least for $K(t)=[h_K+th_C]$ with arbitrary $C\in \mathcal{K}_{(0)}^n$ if $K$ is quasi-smooth.

\begin{lemma}
\label{Wulffvariation}
Let $K,C\in \mathcal{K}_{(0)}^n$ be  such that $rC\subset K$ for some $r>0$.
For $t\in(-r,r)$ and  
$K(t)=[h_K+th_C]$,
\begin{description}
\item[(i)] if $K\subset R\,C$ for $R>0$, then
$|h_{K(t)}(u)-h_K(u)|\leq \frac{R}r \,|t|$ for any $u\in S^{n-1}$ and $t\in(-r,r)$;
\item[(ii)] if $u\in S^{n-1}$ is the exterior normal at some smooth point $z\in\partial K$, then
$$
\lim_{t\to 0}
\frac{h_{K(t)}(u)-h_K(u)}{t}=h_C(u).
$$
\end{description}
\end{lemma}

We will need the condition (\ref{xider2}) in the following rather special setting taken from Bianchi, B\"or\"oczky, Colesanti, Yang \cite{BBCY}.

\begin{lemma}
\label{Wulffvariation-inwards}
Let $K$  be a convex body with $rB^n\subset{\rm int}\, K$ for $r>0$, let $\omega\subset S^{n-1}$ be closed, and if $t\in[0,r)$, then let 
$$
K(t)=[h_K-\mathbf{1}_\omega]=\{x\in K:\,\langle x,u\rangle\leq h_{K}(u)-t\mbox{ \ \ for every $u\in\omega$}\}.
$$
\begin{description}
\item[(i)] We have $\lim_{t\to 0^+}
\frac{h_{K(t)}(u)-h_K(u)}{t}$ exists and is non-positive for all $u\in S^{n-1}$, and  if $u\in\omega$, then even 
$\lim_{t\to 0^+} \frac{h_{K(t)}(u)-h_K(u)}{t}\leq -1$.
\item[(ii)] If $S_K(\omega)=0$, then $\lim_{t\to 0^+}
\frac{V(K(t))-V(K)}{t}=0$.
\end{description}
\end{lemma}


\begin{prop}
\label{Kequasismooth}
For $\varepsilon\in(0,\delta)$, $K^\varepsilon$ is quasi-smooth.
\end{prop}
\proof We suppose that $K^\varepsilon$ is not quasi-smooth, and seek a contradiction.
It follows that $\mathcal{H}^{n-1}(X)>0$ for $X=S^{n-1}\backslash\nu_{K^\varepsilon}(\partial' K^\varepsilon)$, therefore there exists closed $\omega\subset X$ such that $\mathcal{H}^{n-1}(\omega)>0$. Since
$\nu_{K^\varepsilon}^{-1}(\omega)\subset \partial K^\varepsilon\backslash\partial' K^\varepsilon$,
we deduce that $S_{K^\varepsilon}(\omega)=0$.

We may assume that $\xi(K^\varepsilon)=o$ and $rB^n\subset K^\varepsilon\subset RB^n$ for $R>r>0$.
As in Lemma~\ref{Wulffvariation-inwards}, if $t\in[0,r)$, then we define
$$
K(t)=[h_{K^\varepsilon}-\mathbf{1}_\omega]=\{x\in K^\varepsilon:\,\langle x,u\rangle\leq h_{K}(u)-t\mbox{ \ \ for every $u\in\omega$}\}.
$$
Clearly, $K(0)$ equals $K^\varepsilon$. We define $\alpha(t)=V(K(t))^{-1/n}$, and hence $\alpha(0)=1$, and  
Lemma~\ref{Wulffvariation-inwards} (ii) yields that $\alpha'(0)=0$.

We set $\widetilde{K}(t)=\alpha(t)K(t)$, and hence
$\widetilde{K}(0)=K^\varepsilon$ and $V(\widetilde{K}(t))=1$ for all $t\in [0,r)$.
In addition, we consider 
$h(t,u)=h_{K(t)}(u)$ and $\tilde{h}(t,u)=h_{\widetilde{K}(t)}(u)=\alpha(t)h(t,u)$
for $u\in S^{n-1}$ and $t\in[0,r)$.
Since $[h_{K^\varepsilon}-th_{B^n}]\subset K(t)$, Lemma~\ref{Wulffvariation} (i) yields that
$|h(t,u)-h(0,u)|\leq \frac{R}r\,t$ for $u\in S^{n-1}$ and $t\in[0,r)$.
Hence $\alpha'(0)=0$ implies that there exist $c>0$ and $t_0\in(0,r)$ such that
$|\tilde{h}(t,u)-\tilde{h}(0,u)|\leq c\,t$ for $u\in S^{n-1}$ and $t\in[0,t_0)$.
Applying $\alpha(0)=1$, $\alpha'(0)=0$ and Lemma~\ref{Wulffvariation-inwards} (i), we deduce that
\begin{eqnarray*}
\partial_1 \tilde{h}(0,u)&=&\lim_{t\to 0^+}\frac{\tilde{h}(t,u)-\tilde{h}(0,u)}t=
\lim_{t\to 0^+}\frac{h(t,u)-h(0,u)}t\leq 0\mbox{ \ exists for all $u\in S^{n-1}$},\\
\partial_1 \tilde{h}(0,u)&\leq &-1\mbox{ \ for all $u\in \omega$}
\end{eqnarray*}
As $\psi_\varepsilon$ is positive and monotone decreasing, $f>\tau_1$ and  $\mathcal{H}^{n-1}(\omega)>0$,
 Corollary~\ref{center-irrelevant} implies that
\begin{eqnarray*}
\left.\frac{d}{dt}\Phi_\varepsilon(\widetilde{K}(t),\xi(\widetilde{K}(t)))\right|_{t=0^+}&=&
-\int_{S^{n-1}}\partial_1\tilde{h}(0,u)\cdot \psi_\varepsilon(h_K(u))\,f(u)\,d\mathcal{H}^{n-1}(u)\\
&\geq&
-\int_\omega (-1)\psi_\varepsilon(R)\tau_1\,d\mathcal{H}^{n-1}(u)>0.
\end{eqnarray*}
Therefore $\Phi_\varepsilon(\widetilde{K}(t),\xi(\widetilde{K}(t)))>
\Phi_\varepsilon(K^\varepsilon,\xi(K^\varepsilon))$ for small $t>0$. This contradicts the definition of $K^\varepsilon$ 
and concludes the proof.\hfill Q.E.D.\\

For $\varepsilon\in(0,\delta)$, we define
\begin{equation}
\label{lambdaepsilon}
\lambda_\varepsilon=\frac1n\int_{S^{n-1}}h_{K^\varepsilon-\xi(K^\varepsilon)}\cdot
\psi_\varepsilon(h_{K^\varepsilon-\xi(K^\varepsilon)}) \cdot f\,d\mathcal{H}^{n-1}.
\end{equation}

\begin{prop}
\label{Euler-Lagrange}
For $\varepsilon\in(0,\delta)$, 
$\psi_\varepsilon(h_{K^\varepsilon-\xi(K^\varepsilon)}) \cdot f\,d\mathcal{H}^{n-1}=
\lambda_\varepsilon\,dS_{K^\varepsilon}$ as measures on $S^{n-1}$.
\end{prop}
We omit the proof of this result since it  is very similar to that of \cite[Proposition 6.1]{BBCY}, given by the authors and Yang for the $L_p$ Minkowski problem, with $-\lambda_\varepsilon$, $-\psi_\varepsilon$, Lemma~\ref{Alexandrov}, Lemma~\ref{Wulffvariation}, Corollary~\ref{center-irrelevant},  and \cite{SCH}  replacing respectively $\lambda_\varepsilon$, $\varphi'_\varepsilon$, Lemma 5.2, Lemma 2.3, Corollary 3.6 and [72].

\section{The proof of Theorem~\ref{psi0pn} when $f$ is bounded and bounded away from zero}
\label{secmaintheorembounded}

In this section, again let $0<\tau_1<\tau_2$,  let the real function $f$ on $S^{n-1}$ satisfy
$\tau_1<f<\tau_2$, and let $\varphi$ be the continuous function on $[0,\infty)$ of Theorem~\ref{psi0pn}.
We use the notation developed in Section~\ref{secenergy}, and hence 
$\psi:(0,\infty)\to(0,\infty)$ and $\psi=1/\varphi$.

To ensure that a convex body is "fat" enough in Lemma~\ref{Kebounded} and later, the following
observation is useful:

\begin{lemma}
\label{centroidinradius}
If $K$ is a convex body in $\R^n$ with $V(K)=1$ and $K\subset \sigma(K)+RB^n$
for $R>0$, then
$$
\sigma(K)+r B^n\subset K\mbox{ \ for $r=\frac1{c\kappa_{n-1}}n^{-3/2}\,R^{-(n-1)}$.}
$$
\end{lemma}
\proof Let $z_0+r_0 B^n$ be a largest ball in $K$. According to the Steinhagen
theorem \cite[Theorem 50]{Eggl}, there exists $v\in S^{n-1}$ such that
$$
|\langle x-z_0,v\rangle |\leq c\sqrt{n}r_0 
\mbox{ \ for $x\in K$,}
$$
where $c$ is a positive universal constant. It follows that $1=V(K)\leq c\sqrt{n}r_0 \kappa_{n-1}R^{n-1}$,
thus $r_0 \geq \frac1{c\kappa_{n-1}}n^{-1/2}\,R^{-(n-1)}$.  Since
$\sigma(K)+\frac{r_0}{n}\,B^n\subset K$ by
$-(K-\sigma(K))\subset n(K-\sigma(K))$, we may choose
$r=\frac1{c\kappa_{n-1}}n^{-3/2}\,R^{-(n-1)}$.
\hfill Q.E.D.

We recall (compare (\ref{lambdaepsilon})) that if $\varepsilon\in(0,\delta)$ and $\xi(K^\varepsilon)=o$, then
$\lambda_\varepsilon$ is defined by 
\begin{equation}
\label{lambdaepsilon0}
\lambda_\varepsilon=\frac1n\int_{S^{n-1}}h_{K^\varepsilon}
\psi_\varepsilon(h_{K^\varepsilon})  f\,d\mathcal{H}^{n-1}.
\end{equation}

\begin{lemma}
\label{Kebounded}
There exist $R_0>1$, $r_0>0$ and $\tilde{\lambda}_2>\tilde{\lambda}_1>0$ depending on $f,q,\psi,\aleph$
such that if $\varepsilon\in(0,\delta_0)$ for $\delta_0=\min\{\tilde{\delta},\frac{r_0}2\}$
where $\tilde{\delta}$ comes from (\ref{control}), then 
$\tilde{\lambda}_1\leq \lambda_\varepsilon\leq \tilde{\lambda}_2$ and
$$
\sigma(K^\varepsilon)+r_0B^n\subset K^\varepsilon\subset \sigma(K^\varepsilon)+R_0B^n.
$$
\end{lemma}
\proof According to (\ref{bigpsicondup}), there exists $\aleph_0>0$ depending on $q,\psi,\aleph$ such that if $\varepsilon\in(0,\delta)$ and
$t\in(0,\delta)$, then
$\Psi_\varepsilon(t)\leq\aleph_0 t^q$.
In addition, $\lim_{t\to\infty}\Psi_\varepsilon(t)=0$ by (\ref{bigpsieinfty}),
therefore we may apply 
Lemma~\ref{diameterest}. Since (\ref{control}) provides the condition
$$
\int_{S^{n-1}}\Psi_\varepsilon(h_{K^\varepsilon-\sigma(K^\varepsilon)}) f\,d\mathcal{H}^{n-1}\geq 
\Psi(5)\int_{S^{n-1}}f\,d\mathcal{H}^{n-1}
$$
for any $\varepsilon\in(0,\tilde{\delta})$, we deduce from Lemma~\ref{diameterest} the existence of 
$R_0>0$ such that $K^\varepsilon\subset \sigma(K^\varepsilon)+R_0B^n$ for any
$\varepsilon\in(0,\tilde{\delta})$. In addition, the existence of $r_0$ independent of $\varepsilon$
such that $\sigma(K^\varepsilon)+r_0B^n\subset K^\varepsilon$ follows 
from Lemma~\ref{centroidinradius}.

To estimate $\lambda_\varepsilon$, we assume $\xi(K^\varepsilon)=o$. 
 Let
$w_\varepsilon\in S^{n-1}$ and $\varrho_\varepsilon\geq 0$ be such that 
$\sigma(K^\varepsilon)=\varrho_\varepsilon w_\varepsilon$, and hence  $r_0w_\varepsilon\in K^\varepsilon$.
It follows that $h_{K^\varepsilon}(u)\geq r_0/2$ holds for 
$u\in \Omega(w_\varepsilon,\frac{\pi}3)$, while $K^\varepsilon\subset 2R_0B^n$, $R_0>1$ and the monotonicity of $\psi_\varepsilon$ imply that 
$\psi_\varepsilon(h_{K^\varepsilon}(u))\geq \psi_\varepsilon(2R_0)=\psi(2R_0)$ for all $u\in S^{n-1}$.

We deduce from (\ref{lambdaepsilon0}) that
$$
\lambda_\varepsilon=\frac1n\int_{S^{n-1}}h_{K^\varepsilon}
\psi_\varepsilon(h_{K^\varepsilon})  f\,d\mathcal{H}^{n-1}\geq
\frac1n\cdot\frac{r_0}2\cdot \psi(2R_0)\cdot \tau_1\cdot \mathcal{H}^{n-1}\left(\Omega\left(w_\varepsilon,\frac{\pi}3\right)\right)=\tilde{\lambda}_1.
$$

To have a suitable upper bound on $\lambda_\varepsilon$, we define $\alpha\in(0,\frac{\pi}2)$ with
$\cos\alpha=\frac{r_0}{2R_0}$, and hence
$$
\Omega(-w_\varepsilon,\alpha)=\left\{u\in S^{n-1}:\,
\langle u,w_\varepsilon\rangle\leq \frac{-r_0}{2R_0}\right\}.
$$
A key observation is that if $u\in S^{n-1}\backslash \Omega(-w_\varepsilon,\alpha)$, then 
$\langle u,w_\varepsilon\rangle> -\frac{r_0}{2R_0}$ and $\varrho_\varepsilon\leq R_0$ imply 
$$
h_{K^\varepsilon}(u)\geq \langle u,\varrho w_\varepsilon+r_0u\rangle\geq r_0-\frac{r_0\varrho_\varepsilon}{2R_0}\geq r_0/2,
$$
therefore $\varepsilon<\frac{r_0}2$ yields
\begin{equation}
\label{lambdaepsupperr}
\psi_\varepsilon(h_{K^\varepsilon}(u))\leq \psi_\varepsilon(r_0/2)=\psi(r_0/2).
\end{equation}
Another observation is that $K^\varepsilon\subset 2R_0B^n$ implies
\begin{equation}
\label{lambdaepsupperR}
h_{K^\varepsilon}(u)<2R_0\mbox{ \ for any $u\in S^{n-1}$.}
\end{equation}
It follows directly from (\ref{lambdaepsupperr}) and (\ref{lambdaepsupperR}) that
\begin{equation}
\label{lambdaepsupper-}
\int_{S^{n-1}\backslash \Omega(-w_\varepsilon,\alpha)}
 h_{K^\varepsilon}\psi_\varepsilon(h_{K^\varepsilon})f\,d\mathcal{H}^{n-1}\leq  (2R_0)\psi(r_0/2) \tau_2n\kappa_n.
\end{equation}

However, if $u\in \Omega(-w_\varepsilon,\alpha)$, then 
 $\psi_\varepsilon(h_{K^\varepsilon}(u))$ can be arbitrary large as $\xi(K^\varepsilon)$ can be arbitrary close to $\partial K^\varepsilon$ if $\varepsilon>0$ is small, and 
hence we transfer the problem to the previous case 
$u\in S^{n-1}\backslash \Omega(-w_\varepsilon,\alpha)$ using
Corollary~\ref{intcond}. 
First applying $\langle u,-w_\varepsilon\rangle\geq \frac{r_0}{2R_0}$ for $u\in \Omega(-w_\varepsilon,\alpha)$,
then Corollary~\ref{intcond}, and after that  $\langle u,w_\varepsilon\rangle\leq 1$,
$f\leq \tau_2$
  and (\ref{lambdaepsupperr}) implies
\begin{eqnarray*}
\int_{\Omega(-w_\varepsilon,\alpha)}\psi_\varepsilon(h_{K^\varepsilon}(u))f(u)\,d\mathcal{H}^{n-1}(u)&\leq&
\frac{2R_0}{r_0}\int_{\Omega(-w_\varepsilon,\alpha)}
\langle u,-w_\varepsilon\rangle \psi_\varepsilon(h_{K^\varepsilon}(u))f(u)\,d\mathcal{H}^{n-1}(u)\\
&=& \frac{2R_0}{r_0}\int_{S^{n-1}\backslash \Omega(-w_\varepsilon,\alpha)}
\langle u,w_\varepsilon\rangle \psi_\varepsilon(h_{K^\varepsilon}(u))f(u)\,d\mathcal{H}^{n-1}(u)\\
&\leq & \frac{2R_0}{r_0}\cdot 
\psi\left(\frac{r_0}{2}\right)\tau_2n\kappa_n.
\end{eqnarray*}
Now (\ref{lambdaepsupperR}) yields
$$
\int_{\Omega(-w_\varepsilon,\alpha)}
h_K\psi_\varepsilon(h_{K})f\,d\mathcal{H}^{n-1}\leq \frac{(2R_0)^2}{r_0}\cdot
\psi\left(\frac{r_0}{2}\right)\tau_2n\kappa_n,
$$
which estimate combined with (\ref{lambdaepsupper-}) leads to
$\lambda_\varepsilon< \left(\frac{(2R_0)^2}{r_0}+2R_0\right)\psi(\frac{r_0}{2})\tau_2n\kappa_n$.
In turn, we conclude Lemma~\ref{Kebounded}. Q.E.D. \\

Now we prove Theorem~\ref{psi0pn} if $f$ is bounded and bounded away from zero.

\begin{theorem}
\label{psi0pnfbounded}
For $0<\tau_1<\tau_2$,  let the real function $f$ on $S^{n-1}$ satisfy
$\tau_1<f<\tau_2$, and let $\varphi:[0,\infty)\to [0,\infty)$  be increasing and continuous
satisfying  $\varphi(0)=0$, 
$\liminf_{t\to 0^+}\frac{\varphi(t)}{t^{1-p}}>0$,
and $\int_1^\infty \frac1{\varphi}<\infty$. Let  $\Psi(t)=\int_t^\infty\frac1{\varphi}$.
Then there exist $\lambda>0$ and a $K\in\mathcal{K}^n_0$ with $V(K)=1$ such that
$$
f\,d\mathcal{H}^{n-1}=\lambda \varphi(h_K)\,dS_K,
$$
as measures on $S^{n-1}$, and
\begin{equation}
\label{psi0pnnotC1cond}
\int_{S^{n-1}}\Psi(h_{K-\sigma(K)}) f\,d\mathcal{H}^{n-1}\geq 
\Psi(5)\int_{S^{n-1}} f\,d\mathcal{H}^{n-1}.
\end{equation}
In addition, if $f$ is invariant under a closed subgroup $G$ of $O(n)$, then $K$ can be chosen to be invariant under $G$.
\end{theorem}
\proof
We assume that $\xi(K^\varepsilon)=o$ for all $\varepsilon\in(0,\delta_0)$ where
$\delta_0\in(0,\delta]$ comes from Lemma~\ref{Kebounded}. Using the constant $R_0$ of Lemma~\ref{Kebounded}, 
we have that 
\begin{equation}
\label{lambdaepsupperR0}
\mbox{$K^\varepsilon\subset 2R_0B^n$ and \ }h_{K^\varepsilon}(u)<2R_0\mbox{ \ for any $u\in S^{n-1}$
and $\varepsilon \in(0,\delta_0)$.}
\end{equation}

We consider the continuous increasing  function 
$\varphi_\varepsilon:[0,\infty)]\to[0,\infty)$ defined by $\varphi_\varepsilon(0)=0$ and
$\varphi_\varepsilon(t)=1/\psi_\varepsilon(t)$ for $\varepsilon\in(0,\delta)$. We claim that
\begin{equation}
\label{psietophi}
\mbox{$\varphi_\varepsilon$ tends uniformly to $\varphi$ on $[0,2R_0]$ as $\varepsilon>0$ tends to zero.}
\end{equation}
For any small $\zeta>0$, there exists $\delta_\zeta>0$ such that
$\varphi(t)\leq \zeta/2$ for $t\in[0,\delta_\zeta]$. 
We deduce from
 (\ref{psiepsilon}) that if $\varepsilon>0$ is small, then
$|\varphi_\varepsilon(t)-\varphi(t)|\leq \zeta/2$ for $t\in[\delta_\zeta,2R_0]$. However
$\varphi_\varepsilon$ is monotone increasing, therefore 
$\varphi_\varepsilon(t),\varphi(t)\in[0,\zeta]$ for $t\in[0,\delta_\zeta]$, completing the proof of 
(\ref{psietophi}).
 
For any $\varepsilon\in(0,\delta_0)$, it follows from Lemma~\ref{Euler-Lagrange} that
$\psi_\varepsilon(h_{K^\varepsilon})f\,d\mathcal{H}^{n-1}=
\lambda_\varepsilon\,dS_{K^\varepsilon}$ as measures on $S^{n-1}$.
Integrating $g\varphi_\varepsilon(h_{K^\varepsilon})$ for any continuous real function $g$ on $S^{n-1}$,
we deduce that
\begin{equation}
\label{Euler-Lagrange0}
f\,d\mathcal{H}^{n-1}=\lambda_\varepsilon \varphi_\varepsilon(h_{K^\varepsilon})\,dS_{K^\varepsilon}
\end{equation}
as measures on $S^{n-1}$.

Since $\tilde{\lambda}_1\leq \lambda_\varepsilon\leq \tilde{\lambda}_2$  
for some $\tilde{\lambda}_2>\tilde{\lambda}_1$ independent of $\varepsilon$
according to Lemma~\ref{Kebounded}, (\ref{lambdaepsupperR0}) yields the
existence of $\lambda>0$, $K\in\mathcal{K}^n_0$ with $V(K)=1$ 
and  sequence $\{\varepsilon(m)\}$ tending to $0$ such that
$\lim_{m\to\infty}\lambda_{\varepsilon(m)}=\lambda$ and
$\lim_{m\to\infty}K^{\varepsilon(m)}=K$. As $h_{K^{\varepsilon(m)}}$ tends uniformly to $h_K$
on $S^{n-1}$, we deduce that 
$\lambda_{\varepsilon(m)}\varphi_{\varepsilon(m)}(h_{K^{\varepsilon(m)}})$ tends uniformly to 
$\lambda \varphi(h_K)$ on $S^{n-1}$. In addition, $S_{K^\varepsilon(m)}$ tends weakly to $S_K$, thus
(\ref{Euler-Lagrange0}) yields 
$$
f\,d\mathcal{H}^{n-1}=\lambda \varphi(h_K)\,dS_K.
$$
We note that if $f$ is invariant under a closed subgroup $G$ of $O(n)$, 
then each $K^\varepsilon$ can be  chosen to be invariant under $G$
according to Lemma~\ref{extremal}, therefore  $K$ is invariant under $G$ in this case.

To prove (\ref{psi0pnnotC1cond}), if  $\varepsilon\in (0,\delta_0)$, then
 (\ref{control}) provides the condition
\begin{equation}
\label{psi0pnnotC1cond1}
\int_{S^{n-1}}\Psi_\varepsilon(h_{K^\varepsilon-\sigma(K^\varepsilon)}) f\,d\mathcal{H}^{n-1}\geq 
\Psi(5)\int_{S^{n-1}} f\,d\mathcal{H}^{n-1}.
\end{equation}
Now
Lemma~\ref{Kebounded} yields that there exists $r_0>0$  such that
if $\varepsilon\in (0,\delta_0)$, then
$\sigma(K^\varepsilon)+r_0B^n\subset K^\varepsilon$ where $0<\delta_0\leq\frac{r_0}2$. 
In particular, if $u\in S^{n-1}$, then
$h_{K^\varepsilon-\sigma(K^\varepsilon)}(u)\geq r_0$, and hence
we deduce from (\ref{bigpsiupper}) that
\begin{equation}
\label{psi0pnnotC1cond2}
\Psi_\varepsilon(h_{K^\varepsilon-\sigma(K^\varepsilon)}(u))\leq
\Psi(h_{K^\varepsilon-\sigma(K^\varepsilon)}(u))+\frac{\pi}2.
\end{equation}
Since $K^{\varepsilon(m)}-\sigma(K^{\varepsilon(m)})$ tends to $K-\sigma(K)$,
(\ref{bigpsilimit}) implies that if $u\in S^{n-1}$, then
\begin{equation}
\label{psi0pnnotC1cond3}
\lim_{\varepsilon\to 0^+}\Psi_\varepsilon(h_{K^\varepsilon-\sigma(K^\varepsilon)}(u))=
\Psi(h_{K-\sigma(K)}(u)).
\end{equation}
Combining (\ref{psi0pnnotC1cond1}), (\ref{psi0pnnotC1cond2}) and (\ref{psi0pnnotC1cond3})
with Lebesgue's Dominated Convergence Theorem,
 we conclude
 (\ref{psi0pnnotC1cond}), and in turn Theorem~\ref{psi0pnfbounded}. 
\hfill Q.E.D.\\

\section{The proof of Theorem~\ref{psi0pn}}
\label{secmaintheorem}

Let $-n<p<0$, let $f$ be a non-negative non-trivial function in $L_{\frac{n}{n+p}}(S^{n-1})$,
and let $\varphi:[0,\infty)\to[0,\infty)$ be a
monotone increasing continuous function 
satisfying $\varphi(0)=0$,
\begin{eqnarray}
\label{phicond1final}
\liminf_{t\to 0^+}\frac{\varphi(t)}{t^{1-p}}&>&0\\
\label{phicond2final}
\int_1^\infty\frac{1}{\varphi(t)}\,dt&<&\infty.
\end{eqnarray}
We associate certain functions to $f$ and $\varphi$. 
For any integer $m\geq 2$, we define $f_m$ on ${\mathbb S}^{n-1}$ as follows:
$$
f_m(u)=\left\{
\begin{array}{ll}
m&\mbox{ \ if $f(u)\geq m$},\\[0.5ex]
f(u)&\mbox{ \ if $\frac1m<f(u)< m$},\\[0.5ex]
\frac1m&\mbox{ \ if $f(u)\leq\frac1m$}.
\end{array}
\right. 
$$
In particular, $f_m\leq \tilde{f}$ where the function $\tilde{f}:S^{n-1}\to[0,\infty)$ 
in $L_{\frac{n}{n+p}}(S^{n-1})$, and hence in $L_1(S^{n-1})$, is
defined by
$$
\tilde{f}(u)=\left\{
\begin{array}{ll}
f(u)&\mbox{ \ if $f(u)>1$},\\[0.5ex]
1&\mbox{ \ if $f(u)\leq 1$}.
\end{array}
\right. 
$$

As in Section~\ref{secenergy}, using (\ref{phicond2final}), we define the function
$$
\Psi(t)=\int_t^\infty \frac1{\varphi} \mbox{ \ for $t>0$}.
$$
According to (\ref{phicond1final}),
there exist some $\delta\in(0,1)$ and $\aleph>1$ such that
\begin{equation}
\label{phicondaleph}
\frac1{\varphi(t)}< \aleph t^{p-1} \mbox{ \ for $t\in(0,\delta)$}.
\end{equation}
We deduce from Lemma~\ref{bigpsibehave} that there exists 
$\aleph_0>1$ depending on $\varphi$ such that
\begin{equation}
\label{bigpsicondalephfinal}
\Psi(t)< \aleph_0 t^p \mbox{ \ for $t\in(0,\delta)$}.
\end{equation}

For $m\geq 2$, Theorem~\ref{psi0pnfbounded} yields  a $\lambda_m>0$ and a convex body 
$K_m\in \mathcal{K}_{0}^n$ with $\xi(K_m)=o\in{\rm int}\,K_m$, $V(K_m)=1$ such that
\begin{eqnarray}
\label{fmsolution}
\lambda_m\varphi(h_{K_m})\,dS_{K_m}&=&f_m\,d\mathcal{H}^{n-1} \\
\label{fmlower}
\int_{S^{n-1}}\Psi(h_{K_m-\sigma(K_m)}) f_m\,d\mathcal{H}^{n-1}&\geq& 
\Psi(5)\int_{S^{n-1}} f_m\,d\mathcal{H}^{n-1}.
\end{eqnarray}
 In addition, if $f$ is invariant under a closed subgroup $G$ of $O(n)$, then
$f_m$  is also invariant under $G$, and hence
 $K_m$ can be chosen to be invariant under $G$.

Since $f_m\leq \tilde{f}$, and $f_m$ converges pointwise to $f$, 
 Lebesgue's Dominated Convergence theorem yields the existence of $m_0>2$ such that
if $m>m_0$, then 
\begin{equation}
\label{intfmest}
\frac12\int_{S^{n-1}} f<\int_{S^{n-1}} f_m<2\int_{S^{n-1}} f.
\end{equation}
In particular, (\ref{fmlower}) implies
\begin{equation}
\label{tildeflower}
\int_{S^{n-1}}\Psi(h_{K_m-\sigma(K_m)}) \tilde{f}\,d\mathcal{H}^{n-1}\geq 
\frac{\Psi(5)}2\int_{S^{n-1}} f\,d\mathcal{H}^{n-1}.
\end{equation}

We deduce from $V(K_m)=1$, $\lim_{t\to\infty}\Psi(t)=0$, (\ref{bigpsicondalephfinal}), (\ref{tildeflower}) and 
Lemma~\ref{diameterest} that there exists $R_0>0$ independent of $m$ such that
\begin{equation}
\label{KmR0ext}
K_m\subset \sigma(K_m)+R_0B^n\subset 2R_0B^n \mbox{ \ for all $m>m_0$}.
\end{equation}
Since $V(K_m)=1$, Lemma~\ref{centroidinradius} yields some $r_0>0$ independent of $m$ such that
\begin{equation}
\label{KmR0int}
\sigma(K_m)+r_0B^n\subset K_m\mbox{ \ for all $m> m_0$}.
\end{equation}

To estimate $\lambda_m$ from below, (\ref{KmR0ext}) implies that
$$
\int_{S^{n-1}}\varphi(h_{K_m})\,dS_{K_m}\leq \varphi(2R_0)\mathcal{H}^{n-1}(\partial K_m)\leq
\varphi(2R_0)(2R_0)^{n-1}n\kappa_n,
$$
and hence it follows from (\ref{fmsolution}) and (\ref{intfmest}) the existence of
$\tilde{\lambda}_1>0$ independent of $m$ such that
\begin{equation}
\label{lambdamlower}
\lambda_m=\frac{\int_{S^{n-1}} f_m\,d\mathcal{H}^{n-1}}{\int_{S^{n-1}}\varphi(h_{K_m})\,dS_{K_m}}\geq\tilde{\lambda}_1\mbox{ \ for all $m> m_0$}.
\end{equation}

To have a suitable upper bound on $\lambda_m$ for any $m>m_0$, 
we choose $w_m\in S^{n-1}$ and $\varrho_m\geq 0$ such that
$\sigma(K_m)=\varrho_mw_m$. We set $B^{\#}_m=w_m^\bot\cap {\rm int}\,B^n$ and consider the relative open set
$$
\Xi_m=(\partial K_m)\cap \Big(\varrho_mw_m+r_0B^{\#}_m+(0,\infty)w_m\Big).
$$
If $u$ is an exterior unit normal at an $x\in\Xi_m$ for $m> m_0$, then
$x=(\varrho_m+s)w_m+rv$ for $s>0$, $r\in[0,r_0)$ and $v\in w_m^\bot\cap S^{n-1}$,
and hence $\varrho_mw_m+rv\in K_m$ yields
$$
\langle u, (\varrho_m+s)w_m+rv\rangle=h_{K_m}(u)\geq\langle u, \varrho_mw_m+rv\rangle,
$$
implying that $\langle u,w_m\rangle\geq 0$; or in other words, 
$u\in \Omega(w_m,\frac{\pi}2)$. Since the orthogonal projection of $\Xi_m$
onto $w_m^\bot$ is $B^{\#}_m$ for $m> m_0$, we deduce that
\begin{equation}
\label{SKmOmega}
S_{K_m}\left( \Omega\left(w_m,\frac{\pi}2\right)\right)\geq
\mathcal{H}^{n-1}(\Xi_m)\geq \mathcal{H}^{n-1}(B^{\#}_m)=
r_0^{n-1}\kappa_{n-1}.
\end{equation}
On the other hand, if $u\in \Omega(w_m,\frac{\pi}2)$ for $m> m_0$, then
$\varrho_mw_m+r_0u\in K_m$ yields
\begin{equation}
\label{hKmOmega}
h_{K_m}(u)\geq\langle u, \varrho_mw_m+r_0u\rangle \geq r_0.
\end{equation}
Combining (\ref{intfmest}), (\ref{SKmOmega}) and (\ref{hKmOmega}) implies
\begin{equation}
\label{lambdamupper}
\lambda_m=\frac{\int_{\Omega(w_m,\frac{\pi}2)} f_m\,d\mathcal{H}^{n-1}}
{\int_{\Omega(w_m,\frac{\pi}2)}\varphi(h_{K_m})\,dS_{K_m}}\leq
\frac{2\int_{S^{n-1}} f\,d\mathcal{H}^{n-1}}
{\varphi(r_0)r_0^{n-1}\kappa_{n-1}}=\tilde{\lambda}_2\mbox{ \ for all $m> m_0$}.
\end{equation}

Since $K_m\subset 2R_0B^n$ and $\tilde{\lambda}_1\leq \lambda_m\leq \tilde{\lambda}_2$ 
for $m>m_0$ by (\ref{KmR0ext}), (\ref{lambdamlower}) and (\ref{lambdamupper}), there exists
subsequence $\{K_{m'}\}\subset \{K_m\}$ such that  $K_{m'}$ tends to some convex compact set $K$
and $\lambda_{m'}$ tends to some $\lambda>0$. As $o\in K_{m'}$ and $V(K_{m'})=1$ for all $m'$, we have
$o\in K$ and $V(K)=1$. 

We claim that for any continuous function $g:\,S^{n-1}\to\R$, we have
\begin{equation}
\label{gintegral}
\int_{S^{n-1}}g\lambda\varphi(h_{K})\,dS_{K}=\int_{S^{n-1}} gf\,d\mathcal{H}^{n-1}.
\end{equation}
As $\varphi$ is uniformly continuous on $[0,2R_0]$ and $h_{K_{m'}}$ tends uniformly to $h_K$
on $S^{n-1}$, we deduce that
$\lambda_{m'}\varphi(h_{K_{m'}})$ tends uniformly to  $\lambda\varphi(h_{K})$ on $S^{n-1}$.
Since $S_{K_{m'}}$ tends weakly to $S_K$,  we have
$$
\lim_{m'\to\infty}\int_{S^{n-1}}g\lambda_{m'}\varphi(h_{K_{m'}})\,dS_{K_{m'}}=
\int_{S^{n-1}}g\lambda\varphi(h_{K})\,dS_{K}.
$$
On the other hand, $|gf_m|\leq \tilde{f}\cdot \max_{S^{n-1}}|g|$ for all $m\ge 2$, and 
$gf_m$ tends pointwise to $gf$ as $m$ tends to infinity. Therefore Lebesgue's Dominated Convergence Theorem implies that
$$
\lim_{m\to\infty}\int_{S^{n-1}} g f_{m}\,d\mathcal{H}^{n-1}
=\int_{S^{n-1}} g f\,d\mathcal{H}^{n-1},
$$
which in turn yields (\ref{gintegral}) by (\ref{fmsolution}). In turn, we conclude Theorem~\ref{psi0pn}
by (\ref{gintegral}). \hfill Q.E.D.

\end{document}